\documentclass [11pt]{article}
\newtheorem{theorem}{Theorem}

\begin{document}

\title{\bf A strange example concerning Property $N_{p}$} 
\author{{\bf  Elena Rubei }}
\date{}
\maketitle

\begin{abstract} 
We exhibit an example  of a line bundle $M$ on a smooth complex projective 
variety  $Y$ s.t.
$M$ satisfies Property $N_{p}$ for some $p$, the  $p$-module of a
minimal   resolution of the ideal of the embedding of $Y$ by $M$
is nonzero  and  $ M^{2}$ does not
satisfy Property $ N_{p}$.
\end{abstract}

\def\thefootnote{}
\footnotetext{
{\bf Address:}
Elena Rubei, Dipartimento di Matematica ``U. Dini'', via Morgagni 67/A 
50134 Firenze, Italia

\hspace{0.24cm}{\bf E-mail address:} rubei@math.unifi.it

\hspace{0.24cm}{\bf 2000 Mathematical Subject Classification:} 
14C20, 13D02.


}

\def\thefootnote{\arabic{footnote}}

\setcounter{footnote}{0}

Let  $ M $ be a very ample  line bundle on a smooth
complex projective variety $Y$ and let  $\varphi_{M}: Y \rightarrow
{\bf P}(H^{0}(Y, M)^{\ast})$ be  the map  associated to $M$.
We  recall the definition of Property $N_{p}$ of Green-Lazarsfeld
 (see  \cite{Green1}, \cite{Green2} and also \cite{G-L}, \cite{Green3}):

{\em
let $Y$ be a smooth complex projective variety and let $L$ be a very ample
 line bundle on $Y$ defining an embedding $\varphi_{L}: Y \hookrightarrow
{\bf P}={\bf P}(H^{0}(Y,L)^{\ast })$;
set $S= S(L)= Sym^{\ast}H^{0}(L) $,
the homogeneous coordinate ring of the projective space ${\bf P}$, and
consider
the graded $S$-module  $G=G(L)= \oplus_{n} H^{0}(Y, L^{n})$; let $E_{\ast}$
\[ 0 \longrightarrow E_{n} \longrightarrow E_{n-1} \longrightarrow 
... \longrightarrow E_{0} \longrightarrow G \longrightarrow 0\] 
be a minimal graded free resolution of $G$; 
let $E_p= \oplus_{q \geq 0} B_{p,q} \otimes S(-q) $;
the line bundle $L$ satisfies Property
$N_{p}$ ($p \in {\bf N}$) if and only if

\hspace{1cm} $B_{0,0}= {\bf C}$
 
\hspace{1cm}
 $B_{0,q}= 0$ for $q > 0$

\hspace{1cm} $B_{i,q}= 0$    
for $q \neq i+1$ and  $1 \leq i \leq p $.}

(Thus $L$ satisfies Property $N_{0}$  if and only if  $Y \subset
{\bf P}(H^{0}(L)^{\ast })$ is projectively normal, i.e.  $L$
is normally generated; $L$ satisfies Property $N_{1}$
if and only if $L$ satisfies Property $N_{0}$
and the homogeneous ideal $I$ of $Y \subset
{\bf P}(H^{0}(L)^{\ast })$ is generated by quadrics;
$L$ satisfies Property $N_{2}$  if and only if
$L$ satisfies Property $N_{1}$ and the module of syzygies among
quadratic
generators $Q_{i} \in I$ is spanned by relations of the form
 $\sum L_{i}Q_{i}=0$, where $L_{i}$  are linear polynomials;
and so on.)

In \cite{Ru} we considered the problem to see ``how Property
$N_{p}$ propagates through powers'':
 let  $ M $ be a line
bundle on a smooth complex projective variety $Y$ and suppose 
$M$ satisfies
Property $ N_{p}$; let $s \in {\bf N}$; we wondered
 for which $ k$ the line bundle $ M^{s}$ satisfies
Property $ N_{k} $. 
We proved: {\em
 Let $Y$ be a smooth complex projective variety and $M$  a line bundle
on $Y$;
 if $M$ satisfies Property $ N_{p}$ then $M^{s}$ satisfies Property
$N_{p}$ if $s \geq p $}.

Besides,  despite  there are several theorems
suggesting that there is a ``propagation property'' for Property
$N_{p}$, we observed that surprisingly  it
 is not true that if  $M$ satisfies Property $N_{p} $
then  any power of $M$ satisfies
 Property $N_{p}$; the counterexample we exhibited was  the following:
  the line bundle $M= {\cal O}_{ {\bf P}^{2}}(2)$ satisfies
 Property $ N_{p}$ for every $p$ (see  the paper \cite{O-P} and
 the references  given there), but,
   by Ottaviani-Paoletti's Theorem in \cite{O-P},
for instance   the bundle  $M^{2}={\cal O}_{ {\bf P}^{2} } (4)$
doesn't satisfy Property $N_{10} $.

The example can be considered not convincing  by someone because the
module $E_{10}$
of a minimal  resolution of $G({\cal O}_{{\bf P}^{2}}(2))$ 
on  $S({\cal O}_{{\bf P}^{2}}(2))$ 
is $0$.

The product of projective spaces furnishes the following  more 
convincing example.

\vspace{0.4cm} 

{\bf Example.}
Let $r \in {\bf N}$ and  $r \geq 1$.
      
The line bundle
${\cal O}_{P^r  \times P^1} (1,2)$   satisfies Property $N_p$
$\forall p$ (see \cite{O-P}, Theorem 3.6).
 
 On the contrary
  ${\cal O}_{P^r \times P^1} (2,4)$   does
not satisfy Property $N_{10}$, in fact:

${\cal O}_{P^1 \times P^1} (2,4)$   does
not satisfy Property $N_{10}$
 by the following Gallego-Purnapranja's Theorem 
 
\vspace{0.2cm}

{\tt Theorem (Gallego-Purnapranja \cite{G-P})}.
{\it Let $a,b \geq 2$. 
The line bundle ${\cal O}_{{\bf P}^{1} \times {\bf P}^{1} }(a,b)$ 
 satisfies Property $N_p$ if and only if $p \leq 2a+2b-3 $.}

\vspace{0.2cm}
   
and the fact that
 ${\cal O}_{P^1 \times P^1} (2,4)$   does
not satisfy Property $N_{10}$
implies  ${\cal O}_{P^r \times P^1} (2,4)$   does
not satisfy Property $N_{10}$, by  Remark 
in Section 2 of \cite{Green2} (precisely: 
write ${\bf P}^{r} = {\bf P}(V^{\ast})$ and 
write the module $E_p$ of a minimal resolution of 
$G({\cal O}_{P^r \times P^1}(2,4)) $
as $E_p = \oplus_{q \geq 0} B_{p,q} \otimes S(-q)$, as 
in the definition of Property 
$N_p$; we have that $B_{p,q}$ are representations of $GL(V)$;
 thus,  if they are nonzero for $r=n$ then they are nozero also  for $r > n$)
(analogously to the proof of Prop. 1.8 in \cite{O-P}).
 
 If $ r$ is sufficiently great, we have that the
$10$-th module of a minimal  resolution of $G({\cal O}_{P^r \times
P^1}(1,2) )$ is nonzero
(in fact, if $X$ is a projective variety in ${\bf P}^n$ and $\dim X = m $
with $m< n-1$, then the lenght $l$ 
of a minimal resolution of the ideal ${\cal I}_X$ of $X$
 is $\geq n-m$: let  $t$ be s.t. $H^{m}({\cal O}_{X}(t))
\neq 0$; thus $H^{m+1}({\cal I}_{X}(t))
\neq 0$; then by breaking the twist by $t$ of
a minimal resolution of ${\cal I}_{X}$ into 
short exact sequences, we get $l \geq m-n$).

\vspace{1cm}

{\bf Acknowledgements.} I thank  G. Ottaviani for a  helpful 
discussion.
 

{\small
\begin{thebibliography}{Dilloo dilloo 83}

\bibitem[G-P]{G-P} F.J. Gallego, B.P. Purnaprajna {\em Some results on
   rational surfaces and Fano varities}  preprint math.AG/0001107
\bibitem[Gr1]{Green1} M. Green {\em Koszul cohomology and the geometry of
    projective varieties I} J. Differ. Geom. {\bf 20},
     125-171 (1984)
\bibitem[Gr2]{Green2} M. Green {\em Koszul cohomology and the geometry of
    projective varieties II} J. Differ. Geom. {\bf 20},
      279-289 (1984)
\bibitem[Gr3]{Green3}  M. Green {\em Koszul cohomology and geometry},
     in: M. Cornalba et al. (eds),
     Lectures on Riemann Surfaces, World Scientific Press (1989)
\bibitem[G-L]{G-L}  M. Green, R. Lazarsfeld {\em On the
  projective normality
     of complete linear series on an algebraic curve} Invent. math.
     {\bf 83}, 73-90 (1986)
\bibitem[O-P]{O-P} G. Ottaviani, R. Paoletti
     {\em Syzygies of  Veronese embeddings} preprint math.AG/9811131
     to appear in Compositio Mathematica
\bibitem[Ru]{Ru}
    E. Rubei, {\em A note on Property  $N_p$}, Manuscripta Mathematica
    101 (2000), 449-455.

\end{thebibliography}

}
\end{document}